\documentclass[11pt]{article}
\newcommand{\bc}{\begin{center}}
\newcommand{\ec}{\end{center}}
\newcommand{\be}{\begin{equation}}
\newcommand{\ee}{\end{equation}}
\newcommand{\bea}{\begin{eqnarray}}
\newcommand{\eea}{\end{eqnarray}}
\newcommand{\ba}{\begin{array}}
\newcommand{\ea}{\end{array}}
\newcommand{\edc}{\end{document}}

\def\f{\varphi}

\def\g{\gamma}
\def\w{\omega}
\def\O{\Omega}
\def\S{\Sigma}

\def\b{\beta}
\def\G{\Gamma}

\def\e{\epsilon}
\def\s{\sigma}
\def\m{\mu}

\def\l{\lambda}
\def\L{\Lambda}

\def\d{\partial}

\def\C{{\cal C}}
\def\H{{\cal H}}
\begin{document}
\thispagestyle{empty}
\begin{center}

{\bf AN ANALYSIS OF ISING TYPE MODELS ON CAYLEY TREE BY A CONTOUR ARGUMENT}\\
\vspace{0.4cm}
U.A. Rozikov\footnote{rozikovu@yandex.ru}\\
{\it Institute of Mathematics, 29, F.Hodjaev str., 700143, Tashkent,
 Uzbekistan}\\
\vspace{0.8cm}
{\bf Abstract}
\end{center}
In the paper  the Ising  model with competing $J_1$ and
$J_2$ interactions with spin values $\pm 1$, on a Cayley
tree of order 2 (with 3 neighbors) is considered .
We  study the structure of the
ground states and verify the Peierls condition for the model.
Our second result gives description of Gibbs measures
for ferromagnetic Ising  model with $J_1<0$ and $J_2=0$, using a
contour argument which we also develop in the paper.
By the argument we also study Gibbs measures for a natural generalization of
the Ising model. We discuss some open problems and state
several conjectures.\\[2mm]


\section{Introduction}
\large

The Ising model, with two values of spin $\pm 1$  was
considered in [Pr],[Za] and became actively
researched in the 1990's and afterwards (see for example [BG],
[BRZ], [BRSSZ]).

In the paper we consider an Ising model on a Cayley tree with
competing interactions and some a natural generalization
of the  model. The goal of the paper is to study of ground
states  and Gibbs measures of the model. The method of our
investigation is a contour method on the Cayley tree, which we will
develop here.

Contour methods have been used in the mathematical physics
community for many years. In the simplest application, one
first rewrites the model under consideration in terms of
contour representing the boundaries between regions where
the spin    variable in question is constant, and then uses
a so-called Peierls argument to show that large contours are rare,
thus proving that the leading configurations consist of
large oceans of the one spin value, with only small islands
of minority spins ([M],[S]). The techniques of this method
is globally known as Pirogov-Sinai theory or contour arguments.
 This technique was pioneered by
Peierls [P] in his study of the Ising model, later formalized more
precisely by Griffiths and Dobrushin [GD] . The original argument
benefited from the particular symmetries of the Ising model. The
adaptation of the method to the treatment of non-symmetric models
is not trivial, and was developed by Pirogov and Sinai [M], [S], 
[PS1], [PS2], [Z1]. Later, a
particularly enlightening alternative version of the argument was
put forward by Zahradnik [Z].

Note, that Pirogov-Sinai theory on Cayley tree is not developed.
The method used for the description of Gibbs measures on Cayley
tree is the method of Markov random field theory and recurrent
equations of this theory (See for example [BG], [GR],[GR1],[GR2], 
[MR], [NR], [R1],[R2],
[RS]). But, if we consider non-symmetric models
on Cayley tree, then the description of Gibbs measures by the method 
becomes a
difficult problem: in this situation , a nonlinear operator $W$
that maps $R^r$ (for some $r\geq 1$) into itself appears and the
problem is then to describe the fixed points of this operator.
Also implementing this method it is very difficult to prove
extremity of a disordered Gibbs measure.  This problem is not
easy even for symmetric models on Cayley tree, which have been
studied in [BRZ] for Ising model and in [GR1]
for Potts model on Cayley tree.

Note, that extremal Gibbs measures are important , since they
describe the possible macro states of physical system. The Gibbs
measures of models on $Z^d$ described using Pirogov -Sinai theory
are automatically extremal. So it is crucial to develop
Pirogov-Sinai theory on Cayley tree.

\section{Definitions and preliminary results}

The Cayley tree $\Gamma^k$ (See [Ba]) of order $ k\geq 1
$ is an infinite tree, i.e., a graph without cycles, from each
vertex of which exactly $ k+1 $ edges issue. Let $\Gamma^k=(V, L,
i)$ , where $V$ is the set of vertexes of $ \Gamma^k$, $L$ is the
set of edges of $ \Gamma^k$ and $i$ is the incidence function
associating each edge $l\in L$ with its endpoints $x,y\in V$. If
$i(l)=\{x,y\}$, then $x$ and $y$ are called {\it nearest
neighboring vertexes}, and we write $l=<x,y>$.
 The distance $d(x,y), x,y\in V$
on the Cayley tree is defined by the formula
\bea
d(x,y)=\min\{ d | \exists x=x_0,x_1,...,x_{d-1},x_d=y\in V
\ \mbox{such that} \nonumber \\
\mbox{the pairs} \ <x_0,x_1>,...,<x_{d-1},x_d>
\ \mbox{are nearest neighboring vertexes} \}.\nonumber
\eea

For the fixed $x^0\in V$ we set

$$ W_n=\{x\in V\ \ |\ \  d(x,x^0)=n\}, $$
$$ V_n=\cup_{m=1}^n W_m=\{x\in V\ \ | \ \  d(x,x^0)\leq n\}, $$
$$ L_n=\{l=<x,y>\in L \ \ |\ \  x,y\in V_n\}. $$
 Denote $|x|=d(x,x^0)$, $x\in V$.

A collection of the pairs $<x,x_1>,...,<x_{d-1},y>$ is called a
{\sl path} from $x$ to  $y$ and we write $\pi(x,y)$ .
 We write $x<y$ if
the path from $x^0$ to $y$ goes through $x$.

It is known that there exists a one-to-one
correspondence between the set  $V$ of vertexes of the Cayley tree
of order $k\geq 1$ and the group $G_{k}$ of the free products of
$k+1$ cyclic  groups  of the second order with generators
$a_1, a_2,..., a_{k+1}$.

Let us define a group structure on the group $\G_{k}$ as follows.
Vertices which corresponds to the "words" $g,h\in G_{k}$ are
called nearest neighbors and are connected by an edge if either
$g=ha_i$ or $h=ga_j$ for some $i$ or $j$. The graph thus defined
is a Cayley  tree of order $k$.

Consider a   left (resp. right) transformation shift on $G_{k}$
defined as: for $ g_0\in G_{k}$  we  put
$$
T_{g_0}h=g_0h \ \ (\textrm{resp.}\ \  T_{g_0}h=hg_0,) \ \  \forall
h\in G_{k}.
$$
It is easy  to see that  the set  of all left  (resp. right) shifts on
$G_{k}$ is  isomorphic to the group $G_{k}$.

\subsection{Configuration space and the model}

We consider models where the spin takes values in the set
$\Phi=\{-1,1\}$ . A {\it configuration} $\s$ on $V$ is then
defined as a function
 $x\in V\to\s(x)\in\Phi$; the set of all configurations coincides with
$\Omega=\Phi^{V}$. Assume on $\O$ the group of spatial shifts acts.
We define  a {\it periodic configuration} as a configuration $\s(x)$
which is invariant under a subgroup of shifts $G^*_k\subset G_k$
of finite index. For a given periodic configuration  the index of 
the subgroup is called the {\it period of the configuration}. A configuration
 that is invariant with respect to all
shifts is called {\it translational-invariant}.

The Hamiltonian  of the Ising model with competing
interactions has the form
$$
H(\s)=J_1 \sum\limits_{<x,y>}\s(x)\s(y)+
J_2 \sum\limits_{x,y\in V: \ \ d(x,y)=2}{\s(x)\s(y)}  \eqno (1)
$$
where $J_1, J_2\in R$ are coupling constants and
$\s\in \Omega$.

\subsection{ Gibbs measure}

We consider a standard $\s$-algebra ${\cal F}$ of subsets of
$\Omega$ generated by cylinder subsets, all probability measures
are considered on $(\Omega,{\cal F})$. A probability measure $\m$
is called a {\it a Gibbs measure} (with Hamiltonian $H$) if it
satisfies the DLR equation: $\forall n=1,2,...$ and $\s_n\in
\Phi^{V_n}$: $$ \m\bigg(\{\s\in \O :
\s|_{V_n}=\s_n\}\bigg)=\int_{\O}\m(d\w)\nu^{V_n}_{\w
|_{W_{n+1}}}(\s_n) $$
where $\nu^{V_n}_{\w |_{W_{n+1}}}$ is the
conditional probability $$ \nu^{V_n}_{\w
|_{W_{n+1}}}(\s_n)=Z^{-1}(\w |_{W_{n+1}})\exp(-\b
H(\s_n||\w|_{W_{n+1}})). $$ where $\b>0$. Here $\s_n|_{V_n}$ and
$\w|_{W_{n+1}}$ denote the restriction of $\s,\w\in \O$ to $V_n$
and $W_{n+1}$ respectively. Next, $\s_n$ is a
configuration in $V_n$ and $ H(\s_n||\w|_{W_{n+1}})$ is defined as
the sum $H(\s_n)+U(\s_n,\w|_{W_{n+1}})$ where $$ H(\s_n)=J_1
\sum\limits_{<x,y>:\ \ x,y\in V_n}{\s(x)\s(y)} +J_2
\sum\limits_{d(x,y)=2:\ \  x,y\in V_n}{\s(x)\s(y)} \eqno(2)$$ $$
U(\s_n,\w|_{W_{n+1}})=J_1\sum\limits_{<x,y>:\ \  x\in V_n,
\ \ y\in W_{n+1}}{\s(x)\w(y)}+ $$
$$ +J_2 \sum\limits_{d(x,y)=2:\ \  x\in
V_n,\ \  y\in W_{n+1}}{\s(x)\w(y)}
+J_2 \sum\limits_{d(x,y)=2:\ \  x, y\in W_{n+1}}{\w(x)\w(y)} $$
Finally, $Z(\w |_{W_{n+1}})$
stands for the partition function in $V_n$ with the boundary
condition $\w|_{W_{n+1}}$: $$ Z(\w
|_{W_{n+1}})=\sum_{\tilde\s_n\in\Phi^{V_n}}\exp(-\b
H(\tilde\s_n||\w|_{W_{n+1}}). $$ It is known (see [S])
that for any sequence $\w^{(n)}\in \O$, any limiting point of the
measures $\tilde\nu^{V_n}_{\w^{(n)}|_{W_{n+1}}}$ is a Gibbs
measure. Here $\tilde\nu^{V_n}_{\w^{(n)}|_{W_{n+1}}}$ is a measure
on $\O$ such that $\forall n'>n$: $$
\tilde\nu^{V_n}_{\w^{(n)}|_{W_{n+1}}}\bigg(\{\s\in\O :
\s|_{V_{n'}}=\s_{n'}\}\bigg)= \left\{ \ba{ll}
\nu^{V_n}_{\w^{(n)}|_{W_{n+1}}}(\s_{n'}|_{V_n}), \ \ \mbox{if} \ \
\s_{n'}|_{V_{n'}\setminus V_n}=\w^{(n)}|_{V_{n'}\setminus
V_n}\\[2mm] 0, \ \  \ \ \mbox{otherwise}. \ea \right. $$

\section{Ground states}

In the sequel for the simplicity we will consider Cayley tree of order two
i.e k=2.

The ground states for models on the cubic lattice $Z^d$ it was studied in
[GD], [K], [HS], [PS1], [PS2].

For a pair of configurations $\s$ and $\f$ that coincide almost
everywhere, i.e. everywhere except for a finite number of positions
, we consider a relative Hamiltonian $H(\s,\f)$, the difference between
the energies of the configurations $\s, \f$ of the form

$$
H(\s,\f)=J_1\sum_{<x,y>}(\sigma(x)\sigma(y)-\f(x)\f(y))+
 J_2\sum_{x,y\in V: d(x,y)=2}(\sigma(x)\sigma(y)-\f(x)\f(y)), \eqno (3)
$$
where $J=(J_1,J_2)\in R^2$ is an arbitrary fixed parameter.

Let $M$ be the set of unit balls $V_1$ with vertexes in $V$.
We call the restriction of the configuration $\s$ to the ball
$b\in M$ the {\it bounded configuration} $\s_b .$ We shall say that
two boundary configurations $\s_b$ and $\s'_{b'}$ belong to the same
class if one of them can be obtained from the other by replacing all
 values of the configuration $\s'_{b'}$ or $\s_b$ by the opposite values
and up to  any motion in $G_2$.

It is easy to show that the number of such classes is 4. Let
 $b=\{e, a_1,a_2,a_3\}$ be the ball with the center $e$, where $e$ is
the identity of $G_2$ and $a_i, i=1,2,3$ are the generators of the group.
Then $\s^1_b=\{-1,-1,-1,-1\},\ \ \s^2_b=\{-1,-1,-1,1\},\ \ 
\s^3_b=\{1,-1,-1,-1\},\ \ \s^4_b=\{1,1,-1,-1\}, $ 
are  representatives of each  4 classes.

 By $\C_i$, $ i=1, 2, 3, 4$ we denote the $i-th$ class of configurations.

Define the energy of ball $b$ for configuration $\s$ by
$$
U(\s_b)\equiv U(\s_b, J)={1\over 2}J_1 \sum\limits_{<x,y>,\ \ x,y\in b}
\s(x)\s(y)+
J_2 \sum\limits_{x,y\in b:\ \  d(x,y)=2}{\s(x)\s(y)},  \eqno (4)
$$
where $J=(J_1, J_2)\in R^2.$ 

The function $U(\s_b)$ on every class $\C_i$ takes one and the same
value, i.e., if $\s_b\in \C_i$ and $\s'_{b'}\in \C_i$ then
$U(\s_b)=U(\s'_{b'}),$ and therefore we denote by $U_i$ the value of
$U(\s_b)$ on class $i=1,2,3,4.$ It is easy to see that
$$ U_1={3\over 2}J_1+3J_2,\ \  U_2={1\over 2}J_1-J_2, \ \ 
U_3={-3\over 2}J_1+3J_2,\ \ 
U_4={-1\over 2}J_1-J_2. \eqno(5)$$

\vskip 0.5 truecm
{\bf Lemma 3.1.} {\it The relative Hamiltonian (3) has the form}
$$H(\s, \f)=\sum_{b\in M}(U(\s_b)-U(\f_b)). \eqno (6)$$
\vskip 0.5 truecm
{\bf Proof.} Note that for any two points $x$ and $y$ such that
$<x,y>$ there are exactly 2 unit balls for which $x$ and $y$
are vertexes. Also, for any two points $u$ and $v$ such that
$d(u,v)=2$ there exist a unique ball for which $u$ and $v$
 are vertexes. This completes the proof.

\vskip 0.5 truecm
{\bf Lemma 3.2.} {\it  For any class $\C_i$ and for any bounded configuration
$\s_b\in \C_i$ there exists a periodic configuration $\f$ with
period non exceeding 2 such that $\f_{b'}\in \C_i$ for any
$b'\in M$ and $\f_b=\s_b.$}

\vskip 0.5 truecm
{\bf Proof.} Consider 4  separate cases.

 {\it Case $\C_1.$} In this case configuration $\f$
coincides with translational-invariant configuration $\f^+=\{\f(x)\equiv +1\}$
or $\f^-=\{\f(x)\equiv -1\}$. Thus the period of $\f$ is 1.

{\it Case $\C_2.$} We continue the bounded configuration $\s_b\in \C_2$ to
whole  lattice $\G^2$ by means of shifts through $g\in \H_1=\{x\in G_2:
\ \ \mbox{number of the}\ \  a_1  \ \ \mbox{in} \ \ x \ \ $ $
\mbox{ is even}\}$. Note (see [GR2]) that $\H_1$ is normal subgroup of $G_2$ of
index 2. So we obtain a periodic configuration with period 2
(=index of the subgroup) which we denote by $\f$; then $\f_b=\s_b.$
It is easy to see that all restrictions $\f_b, b\in M$ configuration
$\f$ belong to $\C_2.$  Consequently, $U(\f_b)=U(\f_{b'})$ for any
$b, b'\in M.$

{\it Case $\C_3.$} In this case let us  consider  $ \H_2=\{x\in G_2:
\ \ |x|\ \  \mbox {is even}\}$.   $\H_2$ is also normal
subgroup of $G_2$ of
index 2. The proof is then similar to proof of the case $\C_2$.

{\it Case $\C_4.$} The proof is similar to proof of the case $\C_2$
here one can take   $ \H_{12}=\{x\in G_2:
\ \ \w_1(x)+\w_2(x)\ \  \mbox {is even}\},$ where $\w_i(x)$ is the number of
$a_i$ in $x\in G_2$. Note (see [GR2]) that $\H_{12}$ is also normal
subgroup of $G_2$ of index 2. The lemma is proved.

\vskip 0.5 truecm
{\bf Definition 3.3.} A periodic configuration $a$ is called a {\it ground
state} for the relative Hamiltonian $H$ if
$$ U(a_b)=\min\{U_1, U_2, U_3, U_4\},\ \  \mbox{ for any}
\ \  b\in M. \eqno (7)$$

\vskip 0.5 truecm
{\bf Remark 1.} A ground state can be defined differently
as a periodic configuration such that for any configuration $\s$
that coincides with $a$ almost everywhere $H(a,\s)\leq 0.$
It is easy to see that from definition 3.3 follows second
definition i.e. $H(a,\s)\leq 0.$ In [PS1],[PS2] it was proved
that these two definitions are equivalent for Hamiltonians on $Z^d$.
But there is a problem to prove of the equivalentness for Hamiltonians
on Cayley tree: normally the ratio of the number of boundary sites to the
number of interior sites of a lattices becomes small in the thermodynamic
limit of a large system. For the Cayley tree it does not, since both
numbers grow exponentially like $k^n$.\\

Correspondingly, we make a 
\vskip 0.5 truecm
{\bf Conjecture 1.} The conditions (7) and $H(a, \s)\leq 0$ are equivalent.

\section {A separation of the set of parameter $J\in R^2.$}

We set
$$U_i(J)=U(\s_b,J), \ \ \mbox{if} \ \ \s_b\in \C_i, i=1,2,3,4.$$

The quantity $U_i(J)$ is a linear function of the parameter $J\in R^2.
$ For every fixed $m=1,2,3,4$ we denote by $A_m$ the set of points $J$
such that

$$U_m(J)=\min\{U_1(J), U_2(J), U_3(J), U_4(J)\}. \eqno (8)$$

It is easy to check that
$$A_1=\{J\in R^2: J_1\leq 0;\ \  J_1+4J_2\leq 0\}; $$
$$A_2=\{J\in R^2: J_1\leq 0;\ \  J_1+4J_2\geq 0\}; $$
$$A_3=\{J\in R^2: J_1\geq 0;\ \  J_1-4J_2\leq 0\}; $$
$$A_4=\{J\in R^2: J_1\geq 0;\ \  J_1-4J_2\geq 0\}; $$
and $R^2=\cup_{i=1}^4A_i .$

\section{Description  of ground states}

For every point $J\in R^2$ we divide the sets $A_i, i=1,2,3,4 $
into two classes of sets $A_{m_1}, ...,A_{m_r}$ and $A_{m_{r+1}},...,
A_{m_4}$ such that
$$J\in \cap_{q=1}^rA_{m_q}\ \ \mbox{and} \ \
J\notin \cup_{q=r+1}^4A_{m_q}, \eqno (9)$$
where $1\leq r\leq 4$; $m_q\in \{1, 2, 3, 4\}$ and $m_q\ne m_p$ if
$p\ne q$; $p,q\in \{1, 2, 3, 4\}.$

\vskip 0.5 truecm
{\bf Lemma 5.1.} {\it A periodic configuration $a$ is a ground state for the
relative Hamiltonian $H$ if and only if for $b\in M$ it is true that $a_b\in \C(m_1,...
m_r)$, where $\C(m_1,...,m_r)=\cup^r_{q=1}\C_{m_q}.$}

\vskip 0.5 truecm
{\bf Proof.} The condition (9) corresponds exactly to
$$U_{m_1}=...=U_{m_r}<\min\{U_{m_q}, q=r+1,...,4\}.$$
What we require follows from definition 3.3.

Suppose two unit balls $b$ and $b'$ are neighbors, i.e.,
have a common edge. We shall then say that the two bounded
configurations $\s_b$ and $\s'_{b'}$ are {\it compatible} if they
coincide on the common edge of the balls $b$ and $b'$.

\vskip 0.5 truecm
{\bf Lemma 5.2.} {\it Suppose condition (9) is satisfied and that
for any $\s_b\in \C(m_1,...,m_r)$ and for any ball $b'$
that neighbors $b$ there exists exactly one bounded configuration
$\s'_{b'}$ that belongs to $\C(m_1,...,m_r)$ and is compatible with
$\s_b$; then the period of the ground state $a$ with $a_b\in \C(m_1,...,m_r)$
for any $b\in M$ does not exceed 2.}

\vskip 0.5 truecm
{\bf Proof.} We fix some $b_0\in M.$ By lemma 3.2 there exists a periodic
configuration $\s$ with period not exceed 2 such that $\s_{b_0}=a_{b_0}$
and $\s_b\in \C(m_1,...,m_r)$ for any $b\in M.$ We show that $a=\s.$
If this is not so, there exist two neighboring balls $b$ and $b'$
such that $a_b=\s_b$ and $a_{b'}\ne \s_{b'}$, i.e. for the bounded
 configuration  $a_b\in \C(m_1,...,m_r)$ and for the ball $b'$ that
is the neighbor of $b$ there exist two different bounded configurations
$\s_{b'}$ and $a_{b'}$ belonging to $\C(m_1,...,m_r).$ We have
obtained an assertion that contradicts the condition of lemma 5.2.
Therefore, the period of $a$ does not exceed 2.

\vskip 0.5 truecm
{\bf Lemma 5.3.}{\it Suppose condition (9) is satisfied and that
for some $\s_b\in \C(m_1,...,m_r)$ and for some ball $b'$
that is a neighbor of  $b$ there exists two bounded configurations
$\s^1_{b'}$ and $\s^2_{b'}$ that belongs to $\C(m_1,...,m_r)$ and
 compatible with $\s_b$; then the relative Hamiltonian $H$ has
infinitely many  ground states.}

\vskip 0.5 truecm
{\bf Proof.} In accordance with lemma 3.2, there exist periodic
configurations $\s^1$ and $\s^2$ with period not exceeding 2 whose
restriction to ball $b'$ are  bounded configurations $\s^1_{b'}$
and $\s^2_{b'}$ respectively. We construct the following configuration
$a^t:$
$$a^t(x)=\left\{\begin{array}{ll}
\s^1(x)& \textrm{if $|x|\in [2tp-2t; 2tp)$}\\
\s^2(x)& \textrm{if $|x|\in [2tp; 2tp+2t)$}\\
\end{array}\right.$$
where $t, p\in \{1,2,...\}.$
For any $b\in M$ it is true that $a^t_b\in \C(m_1,...,m_r)$,
and therefore $a^t$ is a ground state. Since $t\in \{1,2,...\}$ is
arbitrary, the number of ground states is infinite. The lemma is proved.

Summaries, we have

\vskip 0.5 truecm

{\bf Theorem 5.4.} {\it For any $J\in R^2\setminus (\{J\in R^2: J_1=0\}
\cup \{J\in R^2: J_1=\pm 4J_2,\ \  J_2\geq 0\})$ the period of a ground
state for the relative Hamiltonian $H$ does not exceed 2. On 
$\{J\in R^2: J_1=0\}\cup \{J\in R^2: J_1=\pm 4J_2,\ \  J_2\geq 0\}$
there are infinitely many ground states.}

\section{The Peierls condition}

By theorem 5.4. it is obvious that there exist not more than
$2^4=16$ ground states with period not exceeding 2.
For every point $x\in V$ we denote by $V_2(x)$ the ball
$$V_2(x)=\{y\in V: d(x,y)\leq 2\}$$
We denote the restriction of configuration $\s$ to $V_2(x)$
by $pr(\s, V_2(x)).$

\vskip 0.5 truecm
{\bf Definition 6.1.} Let $\s^1,...,\s^q$ be the complete set of
all ground states of the relative Hamiltonian $H$ and suppose their period
does not exceed 2. The ball $V_2(x)$ is said to be an improper ball of the
configuration $\s$ if $pr(\s, V_2(x))\ne pr(\s^j, V_2(x))$ for any 
$j=1,...q.$
The union of the improper balls of the configuration $\s$ is called the
{\it boundary of the configuration} and denoted by $\d(\s).$

\vskip 0.5 truecm
{\bf Definition 6.2.} The relative Hamiltonian $H$ with ground states
 $\s^1,...,\s^q$  satisfies the Peierls condition if for any $j=1,...,q$
and any configuration $\s$ coinciding almost everywhere with $\s^j, $
$$H(\s,\s^j)\geq \lambda |\d(\s)|,$$
where $\l$ is a positive constant that does not depend on $\s$,
and $|\d(\s)|$ is the number of unit balls in $\d(\s).$

\vskip 0.5 truecm
{\bf Theorem 6.3.} {\it Suppose the sets $A_m$ are split in such a way
that condition (9) holds. If any bounded configuration
 $\s_b\in\C(m_1,...,m_r)$ and for any $b'$ that is a neighbor
of $b$ there exist exactly one bounded configuration $\s'_{b'}
\in\C(m_1,...,m_r)$ compatible with $\s_b$ then the Peierls
condition is satisfied.}

\vskip 0.5 truecm

{\bf Proof.} It follows from lemma 5.2 that the period of
a ground state not exceed 2, so that there are not more than
$2^4$ ground states. We prove the fulfillment of the Peierls condition.
Suppose $\s$ coincides almost everywhere with the ground state $\s^j$
and $b\in\d(\s).$ Among the four vertexes of the ball $b$ there is a vertex
$x\in b$ such that $V_2(x)\subset \d(\s).$ Indeed, one can take the center
of $b$ for $x$. In the ball $V_2(x)$ there exists a unit ball $b'$
such that $U(\s_{b'})-U(\s^j_{b'})\geq \e$, where $\e=\min\{U_{m_p},
p=r+1,...,4\}-U_{m_1}.$ Since $\s^j$ is a ground state,
$U(\s_b)\geq U(\s^j_b)$ for any $b\in M.$ Thus for any $b\in \d(\s),$
$$\sum_{b'\in \widehat{V}(b)}(U(\s_{b'})-U(\s^j_{b'}))\geq \e,$$
where $\widehat{V}(b)=\cup_{x\in b}V_2(x).$
Since $|\widehat{V}(b)|=21,$
$$H(\s,\s^j)=\sum_{b\in M}(U(\s_b)-U(\s^j_b))=
\sum_{b\in \d(\s)}(U(\s_b)-U(\s^j_b))\geq {\e\over 21}|\d(\s)|.$$
Therefore, the Peierls condition is satisfied for $\lambda={\e\over 21}$.
The theorem is proved.

\vskip 0.5 truecm
{\bf Conjecture 2.} The models is considered here satisfies the
Peierls condition iff the number of ground states is finite 
(cf. [HS], [Pe]).

\vskip 0.5 truecm
{\bf Remark 2.} If one want to prove the conjecture 2 by well known
arguments (see for example  [HS], [K]) then appears the problem 
mentioned in Remark 1. We hope there is an other argument to prove it.

\section{ The existence of two Gibbs measures}

In this section we consider the model (1) with $J_2=0$
and shall prove that there are at least two Gibbs measure
for the model. Note that the result it was proved [Pr] using
theory of Markov random fields and recurrent equations of this theory.
Here we shall use our ``contour method'' on Cayley tree. The existence
of several Gibbs measures (for some values of parameters $J_1, J_2$)
is a mathematical expression of the well-known physical phenomenon-the
coexistence of several aggregate states (or phases) of the matter.

We recall that the Ising model ((1) with $J_2=0$) is a lattice spin
system described by a configuration of ``spins'' inside $\L\subset V$,
$\s=\{\s(x)\in \{-1,1\}, \ \ x\in \L \}, $ in which only neighboring
spins interact. The energy $H_{\L}(\s |\f)$ of the configuration $\s$
in the presence of boundary configuration $\f=\{\f(x), x\in V\setminus \L\}$
is expressed by the formula
$$H_{\L}(\s |\f)=J_1\sum_{<x,y>:\ \ x,y\in \L}\s(x)\s(y)+
J_1\sum_{<x,y>:\ \ x\in \L, \ \ y\in V\setminus\L}\s(x)\f(y). \eqno (10)$$

We consider the case $J_1<0$-the so-called ferromagnetic Ising model.
For simplicity (without lose of the generality) we put $J_1=-1$.
By theorem 5.4 and theorem 6.3 the Hamiltonian (10) satisfies the Peierls
condition with two ground states $\f_+\equiv 1$ and $\f_-\equiv -1.$
The particular interest for the Ising model is due primarily to the fact
that for the Ising model with $k=2$ the
thermodynamic functions can be calculated explicitly.

The Gibbs measure on the space $\O_{\L}=\{-1,1\}^{\L}$ with
boundary condition $\f$ is defined in the usual way.
$$\m_{\L,\b}(\s/\f)\equiv \m^{\f}_{\L,\b}(\s)={\bf Z}^{-1}(\L,\b,\f)
\exp(-\b H_{\L}(\s |\f)),$$
where ${\bf Z}(\L, \b, \f)$ is the normalizing factor (statistical sum).
The main goal of the section is to prove the following
\vskip 0.5 truecm
{\bf Theorem 7.1.} {\it  For all sufficiently large $\b$ there are at least
two Gibbs measure for the two-dimensional (i.e. $k=2$) ferromagnetic
Ising model on Cayley tree.}
\vskip 0.5 truecm
{\bf Proof.} Let us consider a sequence of balls on $\G^2$
$$V_1\subset V_2\subset ... \subset V_n\subset ..., \ \ \cup V_n=V,$$
and two sequences of boundary conditions outside these balls:
$$\f_{n,+}\equiv 1, n=1,2,..., \ \ \f_{n,-}\equiv -1, n=1,2,...$$
Each of two sequences of measures $\{\m^{\f_{n,+}}_{V_n,\b}, n=1,2,...\}$
and  $\{\m^{\f_{n,-}}_{V_n,\b}, n=1,2,...\}$ contains a convergent
subsequence (See [S, Theorem 1.2]).
We denote the corresponding limits by $\m^+_{\b}, \m^-_{\b}$ for the
 first and second sequence respectively. Our purpose is to show for
 a sufficiently large $\b$ these measures are different.
Now we describe a boundary of configuration which is more simple than
it was defined at the section 6.

Consider $V_n$, let $V'_n\subseteq V_n$, $V'_n=\{t\in V_n: \s(t)=-1\}.$
For any $A\subset V$ denote
$$\d A=\{x\in V\setminus A: \ \ \exists y\in A, <x,y>\}.$$
Consider the graph $G^n=(V'_n,L'_n),$ where
$$
L'_n=\{l=<x,y>\in L:\ \  x,y\in V'_n\}.$$
It is clear, that for a fixed $n$ the graph $G^n$ contains a finite
number of connected subgraphs $G^n_j$ i.e.
$$G^n=\{G^n_1,...,G^n_m\}, \ \ G^n_i=(V'_{n,i}, L^{(i)}_n).$$

\vskip 0.5 truecm
{\bf Definition 7.2.} The boundary $\d V'_n$ is called the {\it boundary
of the configuration} $\s(V_n)=\{\s(x), x\in V_n\}$ and denoted by
$\d(\s(V_n)).$ The set $\d V'_{n,i}, i=1,...,m$ is called the {\it contour}
of the boundary $\d(\s(V_n)).$ The set $V_{n,i}'$ is called the
{\it interior} of the contour $\d V'_{n,i}$ and the set $V'_n$
is called the interior of the boundary $\d (\s(V_n)).$

Note that any collection of contours uniquely determines configuration
 $\s$ inside $V_n$ (for a fixed constant configuration $\f$ outside $V_n$).
Indeed, going from a point $x$ of boundary of $\L$ to $x^0\in V$ (where $x^0$
is a point , which corresponds to $e\in G_2$) through the unique
path $\pi(x,x^0)$ we put $+1$ until of the first point of $\d (\s(\L))$
on $\pi(x,x^0)$, crossing the point we put $-1$ until the second point
of $\d(\s(\L))$ on $\pi(x,x^0)$ and crossing this point we put $+1$ and so on.

\vskip 0.5 truecm
{\bf Lemma 7.3.} {\it Let $\g$ be a fixed contour and
$p_+(\g)=\mu^+_{\b}\{\s: \g\subset \d(\s(V_n))\}.$ Then
$$p_+(\g)\leq \exp\{-2\b|\g |\}, $$
where $|\g |$ stands for a number of elements of the set
$\g .$}

\vskip 0.5 truecm 

{\bf Proof.} Denote 
$$|\d(\s(V_n))|=\sum_{\tilde{\g}: \tilde{\g}\subset \d(\s(V_n))}
|\tilde{\g}|.$$
 For any $\s$, which coincides with
$\f_+\equiv 1$ outside of $V_n$ we have
$$H_{V_n}(\s)=H(\s(V_n))+H(\s(V_n)|\f_+(V\setminus V_n))=
1-|V_{n+1}|+2|\d(\s(V_n))|.$$
Indeed, 
$$H_{V_n}(\s)=-\sum_{<x,y>, \ \ \{x,y\}\cap V_n\ne \emptyset }
\s(x)\s(y)=-\S^+-\S^-,$$
where $\S^+  (\S^-)$ is the part of sum taken for such
$\{x,y\},$ that $\s(x)=\s(y)$, $(\s(x)=-\s(y)).$ It is easy to see 
that  $-\S^-=|\d(\s(V_n))|.$ 
The total number of $\{x,y\}$ such that $\{x,y\}\cap V_n\ne \emptyset$
is equal to $|V_n|+|\d V_n|-1=|V_{n+1}|-1.$

Consequently,
$$\S^+=|V_{n+1}|-1-|\S^-|=|V_{n+1}|-1+\S^-.$$
Thus,
$$H_{V_n}(\s)=1-|V_{n+1}|-2\S^-=1-|V_{n+1}|+2|\d(\s(V_n))|.$$

By definition we have (where $\L=V_n$)
$$p_+(\g)={\sum_{\s(\L):\g\subset \d(\s(\L))}\exp\{-\b H_{\L}(\s)\}\over
\sum_{\s(\L)}\exp\{-\b H_{\L}(\s)\}}=$$
$$={\sum_{\s(\L):\g\subset \d(\s(\L))}\exp\{\b |\L|+ \b |\d \L|-2\b |\d
(\s(\L))|-\b\}\over
\sum_{\s(\L)}\exp\{\b |\L|+ \b |\d \L|-2\b |\d
(\s(\L))|-\b\}}=$$
$${\sum_{\s(\L):\g\subset \d(\s(\L))}\exp\{-2\b |\d
(\s(\L))|\}\over
\sum_{\s(\L)}\exp\{-2\b |\d
(\s(\L))|\}}.$$
Denote
$$F_{\g}=\{\s(\L): \g\subset \d(\s(\L))\},$$
$$F_{\g}^-=\{\s(\L): \g\cap \d(\s(\L))=\emptyset\}.$$
Define the map $\chi_{\g};F_{\g}\to F_{\g}^-$ as following:
for $\s(\L)\in F_{\g}$ we destroy the contour $\g$ changing 
the values $\s(x)$ inside of $\g$ to $+1$. The constructed configuration
is $\chi_{\g}(\s(\L))\in F^-_{\g}.$

It is clear that
$$\d (\s(\L))=\d(\chi_{\g}(\s(\L)))\cup \g ,$$
$$|\d (\s(\L))|=|\d(\chi_{\g}(\s(\L)))|+ |\g |.$$
For a given $\g$ the map $\chi_{\g}$ is one -to-one map.

Further, we can write  
$$p_+(\g)={\sum_{\s(\L)\in F_{\g}}\exp\{-2\b |\d\s(\L)|\}\over
\sum_{\s(\L)}\exp\{-2\b |\d\s(\L)|\}}\leq$$
$${\sum_{\s(\L)\in F_{\g}}\exp\{-2\b |\d
(\s(\L))|\}\over
\sum_{\s(\L)\in F^-_{\g}}\exp\{-2\b |\d
(\s(\L))|\}}=$$
$${\sum_{\s(\L)\in F_{\g}}\exp\{-2\b |\d
(\s(\L))|\}\over
\sum_{\s(\L)\in F_{\g}}\exp\{-2\b |\d
(\chi_{\g}(\s(\L)))|\}}=\exp\{-2\b |\g|\}.$$
The lemma is proved.

\vskip 0.5 true cm 
{\bf Lemma 7.4.} {\it For all sufficiently large $\b$, there is a 
constant $C=C(\b)>0, $ such that
$$\m^+_{\b}\{\s(\L): \ \ |\g|>C\ln |\L|\ \ \mbox{for some} \ \ 
\g\subset \d(\s(\L))\}\to 0,$$
as $|\L|\to \infty.$}
\vskip 0.5 truecm 
{\bf Proof.} Denote by $N_t(r)=|\{\g: \ \ t\in \g, |\g|=r\}|$ -
the number of different contours with $t\in \g$ and $|\g|=r.$ 
Note that $N_t(r)\leq  12^{2r-1}.$

Suppose $\b>\ln 12$, then
$$\m^+_{\b}\{\s(\L): \g\subset \d(\s(\L)), \ \ t\in \g, 
|\g|=r\}\leq 12^{2r-1}\cdot e^{-2\b r}.$$
$$\m^+_{\b}\{\s(\L): \g\subset \d(\s(\L)), \ \ t\in \g, 
|\g|>C_1\ln|\L|\}\leq 
$$
$${1\over 12}\sum_{r\geq C_1\ln |\L|}
12^{2r}\cdot e^{-2\b r}\leq {1\over 12}\sum_{r\geq C_1\ln |\L|}
(144\cdot e^{-2\b })^r=$$
$${(144\cdot e^{-2\b })^{C_1\ln|\L|}\over 12(1-144e^{-2\b})}=
{|\L|^{C_1(\ln 144-2\b)}\over 12(1-144e^{-2\b})},$$
where $C_1$ will be defined later.

Thus, we have 
$$\m^+_{\b}\{\s(\L): \exists \g\subset \d(\s(\L)), \ \  
|\g|>C_1\ln|\L|\}\leq 
 {|\L|^{C_1(\ln 144-2\b)+1}\over 12(1-144e^{-2\b})}.$$

The last expression tends to zero if $|\L|\to\infty$ and $C_1>{1\over
2\b -\ln 144}.$ The lemma is proved.

\vskip 0.5 truecm 
{\bf Lemma 7.5.} {\it  If $e\in \L$. Then uniformly by $\L$  
$$\m^+_{\b}\{\s(\L): \s(e)=-1\}\to 0,$$
as $\b\to \infty.$ }
\vskip 0.5 truecm 
{\bf Proof.} If $\s(e)=-1,$ then $e$ is point for interior 
of some contour, we shall write this as  $e\in \mbox{Int} \g.$ Assume $t\in \g$
and $e\in \mbox {Int} \g$ then for any such contour we have 
$|\g|\geq |t|+2.$

Consequently,
$$\m^+_{\b}\{\s(\L): e\in \mbox{Int}\g, \ \ t\in \g, |\g|<C_1\ln|\L|\}\leq 
$$
$${1\over 12}\sum_{r=|t|+2}^{ C_1\ln |\L|}
(144\cdot e^{-2\b })^r\leq {(144e^{-2\b})^{|t|+2}\over 12(1-144e^{-2\b})}.$$
$$\m^+_{\b}\{\s(e)=-1\}\leq \m^+_{\b}\{\s(\L): e\in \mbox{Int}\g, 
\ \  \g \subset \d(\s(\L))\}\leq 
$$
$${1\over 12}\sum_{|t|=1}^{ C_1\ln |\L|}
{(144\cdot e^{-2\b })^{|t|+2}\over 1-144e^{-2\b}}+
\m^+_{\b}\{\s(\L): \exists \g\subset \d(\s(\L)), \mbox{such that}, 
|\g|\geq C_1\ln|\L|\}\leq 
$$
$$ {12^5 e^{-6\b }\over (1-144e^{-2\b})^2}+
{|\L|^{C_1(\ln 144 - 2\b)+1}\over 12(1-144e^{-2\b})}. \eqno (11)$$

For $|\L|\to \infty$ and $\b\to\infty$ from (11) we get
$\m^+_{\b}\{\s(e)=-1\}\to 0$. The lemma is proved.

Let us continue the proof of theorem 7.1. By lemma 7.5 we have
$$ \m^+_{\b}\{\s(e)=-1\}<{1\over 2}. \eqno(12)$$
Using the similar argument one can prove
$$\m^-_{\b}\{\s(e)=+1\}<{1\over 2}.\eqno (13)$$
By (13) we have 
$$\m^-_{\b}\{\s(e)=-1\}>{1\over 2}. \eqno (14)$$
Thus, from (12) and (14) we have $\m^+_{\b}\ne \m^-_{\b}$. 
The theorem 7.1 is proved.
\vskip 0.5 truecm
{\bf Definition 7.6.} A Gibbs measure $\m_{\b},$ for Hamiltonian $\b H$
is called {\it small deviation} of a fixed configuration $\f\in \O$ if 
for any $n>0$
$$\lim_{\b\to \infty}\sup_{x\in V}\m_{\b}(\s: \s(V_n(x))\ne\f(V_n(x)))=0.
\eqno (15)$$
\vskip 0.5 truecm
{\bf Theorem 7.7.} {\it The Gibbs measures $\m^+_{\b}, \ \ \m^-_{\b}$ of
ferromagnetic Ising model on Cayley tree of order 2 are small
deviations of the ground states $\f_+\equiv 1$ and $\f_-\equiv -1$
respectively.}

\vskip 0.5 truecm
{\bf Proof.} 

$$\m^+_{\b}(\s: \s(V_n(x))\ne\f(V_n(x)))\leq
\sum_{A\subseteq V_n(x)}\m^+_{\b}\{\s: \s(A)\equiv -1\}\leq$$
$$\sum_{A\subseteq V_n(x)}|A|\m^+_{\b}\{\s: \s(e)= -1\}=
C(n)\m^+_{\b}\{\s(e)=-1\},\eqno(16)$$
where $C(n)=\sum_{A\subseteq V_n(x)}|A|.$ It is easy to see
that $C(n)$ depends only on $n$. Using lemma 7.5 from (16)
we get (15).
The theorem is proved.

From the theorem 7.7. it follows an additional information about 
the structure of  the ``typical'' configurations for each of the 
constructed measures. Namely, for the $\m^+_{\b}$ $(\m^-_{\b})$ 
almost every configuration $\s$ is such that on a connected set whose
density on the Cayley tree tends to unity as $\b\to\infty$ the 
configuration $\s$ coincides with $\f_+ \ \ (\f_-)$, and all the connected
components of the set $\{x: \s(x)\ne \f_+(x)\}$ ($\{x: \s(x)\ne \f_-(x)\}$)
are finite.

\section{ Some generalizations}

In this section we consider a generalization of previous
results describing the Ising model. We consider a ``spin'' model
with two values of ``spin'' : $v_1$ and $v_2.$ The state of such a
system on the Cayley tree $\G^2$ is determined by the configuration
$\w=\{\w(x), x\in V\}, \ \ \w(x)=v_1$ or $v_2.$ The energy $H$ of the 
configuration $\w$ inside a finite set $\L\subset V$ is
$$H_{\L}(\w_\L|\w_{V\setminus \L})=
\sum_{<x,y>:\ \ x, y \in \L}\l(\w(x),\w(y))+ 
\sum_{<x,y>:\ \ x \in \L, y\in V\setminus \L}\l(\w(x),\w(y)), \eqno(17)$$
here  $\w_{\L}=\{\w(x), \ \ x\in \L\}$ and $\w_{V\setminus \L}
=\{\w(x), \ \ x\in V\setminus \L \}$ are parts of the configuration $\w$
inside the set and outside it. The interaction $\l(v_i, v_j)=\l_{ij},
\ \ i,j=1, 2$ is given by a matrix ${\cal{M}}=(\l_{ij})_{i,j=1, 2}$ of 
second order. If we set $v_1=-1, v_2=1$ and $\l(v_1,v_2)=-v_1v_2$ 
then we get the ferromagnetic Ising model.

Now we consider two outer constant configurations
$$\bar{\w}^{(i)}_{V\setminus \L}=v_i, \ \ i=1, 2.$$

Denote by $H^{(i)}_\L(\w_\L)$ the energy $H_\L(\w_\L | \bar{\w}^{(i)}_
{V\setminus \L})$ corresponding to the configuration $\bar{\w}^{(i)}_
{V\setminus \L}$ and by $P^{(i)}_{\L, \b}, \ \ i=1, 2$ the corresponding
Gibbs measure.

Denote by $\w^{(i)}_\L$ the configuration $\w_\L$ extended by
$\bar{\w}^{(i)}_{V\setminus \L}$ and by $\d(\w^{(i)}_\L)$ its
boundary as was explained in section 7, so that $\d(\w^{(i)}_\L)$ consists 
of points $x$ of the tree such that $\w(x)=v_i$ and there is at least
one $y\in S_1(x)=\{u\in V: d(x,u)=1\}$, such that $\w(x)\ne \w(y).$

\vskip 0.5 truecm 
{\bf Lemma 8.1.} {\it Let $K$ be a connected subgraph of $\G^2$, such that
$|K|=n$, then $|\d K|=n+2.$}
\vskip 0.5 truecm 
{\bf Proof.} We shall use the induction over $n .$ For $n=1$ and 2 
the assertion is trivial. Assume for $n=m$ the lemma is true
i.e from $|K|=m$ follows $|\d K|=m+2.$ We shall prove the assertion for 
$n=m+1$ i.e. for $\tilde{K}=K\cup \{x\}.$  Since $\tilde{K}$ is connected
graph we have $x\in \d K$ and there is unique $y\in S_1(x)$ such that 
$y\in K.$ Thus $\d \tilde{K}=(\d K\setminus \{x\})\cup 
(S_1(x)\setminus \{y\}).$ Consequently,  
$$|\d \tilde{K}|=|\d K|-1+2=m+3.$$ The lemma is proved.
 
\vskip 0.5 truecm 
{\bf Lemma 8.2.} {\it The energy $H^{(i)}_\L(\w_\L), \ \ i=1, 2$ has the form
$$H^{(1)}_\L(\w_\L)=
(\l_{21}+\l_{22}-2\l_{11})|\d (\w_\L)|+3m(\l_{11}-\l_{22})+
\l_{11}(|V_{n+1}|-1);$$
$$H^{(2)}_\L(\w_\L)=(\l_{12}+\l_{11}-2\l_{22})|\d (\w_\L)|+3m(\l_{22}-\l_{11})+
\l_{22}(|V_{n+1}|-1);$$
where $m\equiv m(\w_\L)=|\{\g: \g\subset \d (\w_\L)\}|$ - the number of 
different contours of 
$\d (\w_\L),$ i.e. $\d(\w_\L)=\{\g_1,...,\g_m\}.$}
\vskip 0.5 truecm 
{\bf Proof.} Using lemma 8.1 and well known fact that if $K$ is a connected
graph then number of edges of $K$ equal $|K|-1$ we have
$$H^{(1)}_\L(\w_\L)=
\sum^m_{i=1}\bigg(\sum_{<x,y>: x,y\in {\rm Int}\g_i}\l(\w(x),\w(y))+$$
$$
\sum_{<x,y>: x\in {\rm Int}\g_i, y\in \g_i}\l(\w(x),\w(y))\bigg)+
\sum_{<x,y>: x,y\in V_{n+1}\setminus \cup^m_{i=1}{\rm Int}\g_i}
\l(\w(x),\w(y))=$$
$$\sum^m_{i=1}\bigg(\l_{22}(|\mbox{Int}\g_i|-1)+\l_{21}|\g_i|\bigg)+
\l_{11}\bigg(|V_{n+1}|-1-\sum_{i=1}^m(|\mbox{Int}\g_i|+|\g_i|-1)\bigg)=$$
$$\l_{22}\sum^m_{i=1}(|\g_i|-3)+\l_{21}\sum^m_{i=1}|\g_i|+
\l_{11}\bigg(|V_{n+1}|-1-2\sum^m_{i=1}|\g_i|+3m\bigg)=$$
$$\l_{22}(|\d(\w_\L)|-3m)+\l_{21}|\d (\w_\L)|+\l_{11}(|V_{n+1}|-1-
2|\d(\w_\L)|+3m)=$$
$$(\l_{21}+\l_{22}-2\l_{11})|\d (\w_\L)|+3m(\l_{11}-\l_{22})+
\l_{11}(|V_{n+1}|-1);$$
The proof of second part of the lemma is similar. The lemma is proved.
\vskip 0.5 truecm 
{\bf Lemma 8.3.} {\it Let $\g$ be a fixed contour and $p^{(i)}_\b(\g)=
P^{(i)}_\b\{\w: \g\subset \d(\w_\L)\}.$ Then
$$p^{(1)}_\b(\g)\leq \exp\{-\b(\l_{21}+\l_{22}-2\l_{11})|\g|-
3\b(\l_{11}-\l_{22})\},$$
and a similar inequality holds for $p^{(2)}_\b(\g).$ }
\vskip 0.5 truecm 
{\bf Proof.}  One can use lemma 8.2 and the simple fact that
 $m(\chi_\g(\w_\g))=m(\w_\g)-1,$ where $\chi_\g$ defined in proof of
lemma 7.3. Then the proof of the lemma 8.3 is similar to proof of the 
lemma 7.3. 

Using lemma 8.3 one can prove an analogies of
lemma 7.4 and 7.5  then also
\vskip 0.5 truecm 
{\bf Theorem 8.4.} {\it For all sufficiently large $\b$ there are
at least two Gibbs measure for the two dimensional model (17)
on the Cayley tree of order 2.}

\vskip 0.5 truecm 
{\bf Theorem 8.5.} {\it The Gibbs measure $P^{(i)}_\b$
of the model (17) on Cayley tree of order 2 is small deviation
of the configuration $\w^{(i)}\equiv v_i, i=1, 2.$}
\vskip 0.5 truecm 

{\bf Acknowledgments.} The work supported by NATO Reintegration Grant
: FEL. RIG. 980771. The final part of this work was done within the scheme 
of Mathematical Fellowship at the Abdus Salam International Center for 
Theoretical Physics (ICTP), Trieste, Italy and the author thank ICTP for
providing financial support and all facilities. I thank IMU/CDE-program
for (travel) support. 
\vskip 0.5 truecm 
{\bf References}

1. [Ba] R.J. Baxter, {\it Exactly  Solved Models in Statistical Mechanics},
(Academic Press, London/New York, 1982).

2.[BRZ] P.M. Bleher, J. Ruiz, V.A. Zagrebnov. On the purity of the limiting 
Gibbs state for the Ising model on the
Bethe lattice, {\it Jour. Statist. Phys.} {\bf 79} : 473-482 (1995).

3.[BG] P.M. Bleher and  N.N. Ganikhodjaev, On pure phases of the Ising model
on the Bethe lattice, {\it Theor. Probab. Appl.} {\bf 35}:216-227
(1990) .

4.[BRSSZ] P.M. Bleher, J. Ruiz, R.H.Schonmann, S.Shlosman and V.A. Zagrebnov,
Rigidity of the critical phases on a Cayley tree, {\it Moscow
Math. Journ}. {\bf 3}: 345-362 (2001).

5. [GD] V.M. Gertsik, R.L.Dobrushin, Gibbs states in lattice models with
two-step interaction, {\it Func. Anal. Appl.} {\bf 8}: 201-211 (1974).

6. [GR] N.N.Ganikhodjaev, U.A.Rozikov, Group representation of the Cayley 
forest and some applications, {\it Izvestiya: Math.} {\bf 67}: 17-27 (2003).

7. [GR1] N.N. Ganikhodjaev and U.A. Rozikov, On disordered phase in the 
ferromagnetic Potts model on the Bethe lattice, {\it Osaka Journ. Math.} {\bf
37}: 373-383 (2000). 

8. [GR2] N.N. Ganikhodjaev and U.A. Rozikov, A description of periodic
extremal Gibbs measures of some lattice models on the Cayley tree,
{\it Theor. Math. Phys.} {\bf 111}: 480-486 (1997).

9. [HS] W. Holsztynski, J. Slawny, Peierls condition and the number of
ground states, {\it Commun. Math. Phys.} {\bf 61}: 177-190 (1978).

10. [K] I.A. Kashapov, Structure of ground states in three-dimensional
Ising model with tree-step interaction, {\it Theor. Math. Phys.} {\bf 33}:
912-918 (1977).

11. [M] R.A. Minlos, {\it Introduction to mathematical statistical physics}
(University lecture series, ISSN 1047-3998; v.19,  2000)

12. [MR] F.M. Mukhamedov , U.A. Rozikov, On Gibbs measures of models
with competing ternary and binary interactions and
corresponding von Neumann algebras.  {\it Jour. of Stat.Phys.} {\bf 114}:
825-848 (2004).

13. [NR] Kh.A. Nazarov, U.A. Rozikov, Periodic Gibbs measures for the 
Ising model with competing interactions, {\it Theor. Math. Phys.}
{\bf 135}: 881-888 (2003)

14. [Pe] E.A. Pecherski, The Peierls Condition is not always satisfied,
{\it Select. Math. Sov.} {\bf 3}: 87-92 (1983/84)

15.[P] R. Peierls, On Ising model of ferro magnetism. {\it Proc. Cambridge
Phil. Soc.} {\bf 32}: 477-481 (1936)

16.[Pr] C. Preston, {\it Gibbs states on countable sets}
(Cambridge University Press, London 1974).

17. [PS1] S.A.Pirogov, Ya.G.Sinai, Phase diagrams of classical lattice
systems,I. {\it Theor. Math. Phys.} {\bf 25}: 1185-1192 (1975)

18. [PS2] S.A.Pirogov, Ya.G.Sinai, Phase diagrams of classical lattice
systems,II. {\it Theor. Math. Phys.} {\bf 26}: 39-49 (1976)

19. [R1] U.A. Rozikov, A description of limit Gibbs measures for
$\l$-models on the Bethe lattice, {\it Siberian Math. J.} {\bf
39}: 373-380 (1998).

20. [R2] U.A. Rozikov,  Representation of trees and its application
 {\it Math. Notes.}
{\bf 72}: 516-527 (2002).

21. [R3] U.A. Rozikov, Partition structures of the group representation of the
Cayley tree into cosets by finite-index normal subgroups and their applications
to the description of periodic Gibbs distributions,  {\it Theor. Math. Phys.}
{\bf 112}: 929-933 (1997).

22. [RS] U.A. Rozikov, Yu.M. Suhov, A hard-core model on a Cayley tree:
an example of a loss network, {\it Queueing Syst.} {\bf 46}: 197-212 (2004)

23. [S] Ya.G. Sinai, {\it Theory of phase transitions: Rigorous Results}
(Pergamon, Oxford, 1982).

24. [Z] M. Zahradnik, An alternate version of Pirogov-Sinai theory.
{\it Comm. Math. Phys.} {\bf 93}: 559-581 (1984)

25. [Z1] M. Zahradnik, A short course on the Pirogov-Sinai theory,
{\it Rendiconti Math. Serie VII} {\bf 18}: 411-486 (1998)

26.[Za] S. Zachary, Countable state space Markov random fields and Markov
chains on trees, {\it Ann. Prob.} {\bf 11}: 894-903 (1983).

\end{document}